\documentclass{article}
\usepackage[utf8]{inputenc}
\usepackage{amsmath}
\usepackage{amssymb}
\usepackage{pifont}
\usepackage{url}
\usepackage{algorithm}
\usepackage{algpseudocode}
\usepackage{graphicx}
\usepackage{multirow}
\def\G1{\hbox{$\displaystyle{\mbox{\ding{172}}}$}}
\def\go{\hbox{$\displaystyle{\mbox{\ding{172}}}$}}
\title{High Precision Differentiation Techniques for Data-Driven Solution of Nonlinear PDEs by Physics-Informed Neural Networks}
\author{Marat S. Mukhametzhanov}

\date{\footnotesize Department of Computer Science, Modeling, Electronics and Systems Engineering (DIMES), University of Calabria, 87036 -- Rende (CS), Italy\\
m.mukhametzhanov@dimes.unical.it\\
October 2022}

\begin{document}

\maketitle
\begin{abstract}
    Time-dependent Partial Differential Equations with given initial conditions are considered in this paper. New differentiation techniques of the unknown solution with respect to time variable are proposed. It is shown that the proposed techniques allow to generate accurate higher order derivatives simultaneously for a set of spatial points. The calculated derivatives can then be used for data-driven solution in different ways. An application for Physics Informed Neural Networks by the well-known DeepXDE software solution in Python under Tensorflow background framework has been presented for three real-life PDEs: Burgers', Allen-Cahn and Schrodinger equations. 
\end{abstract}
\section{Introduction}

Time-dependent Partial Differential Equations (PDEs) arise frequently in real-life applications: e.g., diffusion process of liquid flows (see, e.g., \cite{C9SM01119F}), heat distribution in time (see, e.g., \cite{Raissi:et:al.(2018)}), simulations of nonlinear wave dynamics (see, e.g., \cite{WANG2021113938}), groundwater flow dynamics (see, e.g., \cite{groundwater}), quantum dynamics (see, e.g., \cite{griffiths(2004)}), computational mechanics (see, e.g., \cite{LMLG.CMAME-RIKS.2018}), etc. These applications are very important from both theoretical and practical points of view. For instance, groundwater flow simulations can be used to predict hydro-geological risks, which are crucial for infrastructures located in seismic or unstable regions (see, e.g., \cite{Anderson(2015), Kresic(2007)}. High precision efficient simulations and modeling in this case can be used to predict different risks arising in this field. In this case, numerical models can be used to describe fluid dynamics: e.g., diffusion equation or Burgers' equations (see, e.g., \cite{BASDEVANT198623}).

In order to solve difficult nonlinear PDEs, there exist different approaches: e.g., finite element method (FEM, see, e.g., \cite{Korelc&Wriggers(2016)}) or Isogeometric analysis (IGA, see, e.g., \cite{HughesBook,BordasOverview}). Recently, a wide interest to deep learning has led to another approach for solving challenging PDEs: data-driven solution by Phisics-Informed Neural Networks (PINNs, see, e.g., \cite{cuomo(2022),nature_pinns,RAISSI2019686}). 

Neural networks are usually called as universal approximators allowing one to approximate any function with any given accuracy if there is sufficient data. For this reason, they can be used for the approximation of the solution to differential equations if there is enough information or data about the solution (e.g., sufficiently large data set on boundary and initial conditions and sufficiently large data set of collocation points, i.e., the points inside the domain, which are used to maintain the PDE dynamics minimizing the loss function defined through the equation without initial or boundary conditions, see, e.g., \cite{RAISSI2019686}). 

However, in practice the information about the solution of the PDE can be very low: only initial and boundary conditions on a limited (e.g., onedimensional) domain. In this case, too large set of collocation points can lead to a bias of the solution into the dynamics of the PDE without taking into account real values of the solution and, thus, to a small accuracy with respect to the exact solution. Thus, deterministic methods can be used with the PINNs to obtain additional high precision data to train or tune the deep learning models (see, e.g., \cite{RAISSI2019686}).  As it has been shown in \cite{Leonetti&Mukhametzhanov(2022)}, high precision derivatives for PDEs can also improve the performance of the simulation without necessity of elaboration of complex analytical formulae. In this paper, to overcome the above mentioned issues related to lack of data, new differentiation techniques are proposed for the time-dependent PDEs. The proposed schemes allow to generate much more data about the solution of the PDE giving so much more additional information to the neural networks for data-driven solution.

PINNs are used usually to solve two kinds of problems related to differential equations: data-driven solution and data-driven discovery (see, e.g., \cite{RAISSI2019686}). More specifically, there can be considered the following nonlinear Partial Differential Equation (PDE):
\begin{equation}
    U_t = N[U,\lambda],~U = U(t,x),~t\in [t_0, t_{end}],~x\in D,
    \label{eq:pde_general}
\end{equation}
where the data-driven solution consists of approximation of the solution $U(t,x)$ using the data points on the initial and boundary conditions given all the parameters $\lambda$ of the model, while the discovery problem consists of approximation of the model parameters $\lambda$ given the values of $U(t,x)$ at some points. In this paper, the data-driven solution problem is considered. 

Data-driven solution consists of approximation of the unknown solution $U(t,x)$ by neural networks using the (noisy) observations of the function $U(t,x)$. The loss function, usually, is defined by combinations of the given data and collocations points: e.g., as $L(t,x) = ||V(t,x)-U(t,x)||_{\Gamma_u} + ||f_V(t,x)||_{\Gamma_f}$, where $U(t,x)$ and $V(t,x)$ are the real and approximated solutions, respectively, $f_V(t,x) = V_t - N[V,\lambda]$, while $||\cdot ||_{\Gamma}$ is some norm (e.g., Euclidean) over the set $\Gamma$, which consists of the respective values $(t,x)$ used to define the loss function. Collocation points $(t,x)\in [t_0, t_{end}]\times D$ can be chosen as many as needed, since they are used only to guarantee the PDE (\ref{eq:pde_general}) for the approximated solution $V(t,x)$ without taking into account real values of $U(t,x)$. Differentiation of $V(t,x)$ can be performed, e.g., using Automatic Differentiation (for instance, by Tensorflow). Known observations of $U(t,x)$ can be given, frequently, by initial condition and/or boundary conditions (in the latter case, the respective loss function is also defined). Sometimes, the known observations of the solution $U(t,x)$ can be given as a set of (noisy)  measurements without any analytical and/or software code. In this paper, we consider the case when the initial value condition is given by a function $g(x),~x\in D$, which can be also given by computational procedures, e.g., as a ``black box'' without analytical formulae. 

It becomes clear, for the above mentioned issues, that the availability of the high precision data is crucial for the accuracy of the solution to PDEs using PINNs. However, in the case when there is no additional data, except the initial and boundary conditions, it can be difficult to tune the neural network to produce the results with high precision. In that case, the available data is limited only by the values at $t = t_0$ and, probably, by boundary conditions, which cannot give all necessary information about the solution (e.g., its dynamics around the initial time at $t > t_0$), even if some additional information can be obtained using PINNs along with deterministic algorithms for PDEs (see, e.g., \cite{RAISSI2019686}). Derivatives of $U(t,x)$ at the initial time $t_0$ can provide an important information, which can be used to tune better the neural networks: e.g., adding another loss function on the derivatives $L_d(t,x) = \sum_{i=1}^{n} ||\frac{\partial^{i} V}{\partial t^i}(t_0,x) - \frac{\partial^{i} U}{\partial t^i}(t_0,x)||_{\Gamma_x}$, or trying to obtain a higher precision approximation of $U(t,x)$ around $t_0$ (e.g., by the Taylor expansion) to use it for fine-tuning the network. 

However, since the initial condition is given only at the fixed time $t_0$, it is impossible to differentiate it directly, e.g., by Automatic Differentiation, to obtain high precision derivatives. In the papers \cite{ODE_4,ODE_5}, several algorithms for computation of the higher order derivatives with respect to the time variable have been proposed for Ordinary Differential Equations (ODEs). There are two most important drawbacks regarding the algorithms presented there. The first one is related to the using of the right-side function of the ODE for computation of the derivatives. In particular, the function $f(y)$ from the (autonomous) ODE $y'(t) = F(y)$ should be available to be calculated at any concrete fixed numerical value of $y$. However, this is not the case of PDEs, since the right-side function $N[U,\lambda]$ can depend also on the derivatives $U_x$ and $U_{xx}$ of the (unknown) solution $U(t,x)$ and cannot be calculated for a fixed numerical value of $U$. For this reason, the algorithms proposed in \cite{ODE_4,ODE_5} cannot be used directly to calculate the higher order derivatives of $U$ with respect to $t$. 

The second issue in the above mentioned algorithms consists in the using of a software for working with infinite and infinitesimal quantities (e.g., Levi-Civita field \cite{Flynn&Shamseddine(2020),Shamseddine&Berz(1996)}, the Infinity Computer \cite{Sergeyev_patent,Sergeyev_EMS_survey}, Nonstandard analysis \cite{Robinson}, etc.). The Infinity Computer was used in the above mentioned papers for dealing with infinite and infinitesimal quantities. It can be described briefly as the computer system, which is based on the positional numeral system with the infinite radix $\go$ (defined as the number of elements of the set of natural numbers, i.e., $\go = |\mathbb{N}|$, see, e.g., \cite{Sergeyev_EMS_survey}, section 4.1). A number $C$ on the Infinity Computer is represented as follows (see \cite{Sergeyev_EMS_survey}, section 4.3): 
\begin{equation}
    C = c_0\go^{p_0} + c_1\go^{p_1} + ... + c_n \go^{p_n},
    \label{eq:grossnumber}
\end{equation}
where $c_i \neq 0,~i = 0,...,n,$ are finite floating-point numbers called grossdigits, while $p_i,~i=0,...,n,$ are called grosspowers and can be finite, infinite\footnote{It should be noted that only finite grosspowers are used so far at all sources related to numerical computations, since there is not defined a procedure to compare different numbers with infinite and/or infinitesimal grosspowers as it has been shown in \cite{Gutman:et:al.(2017)}, Fact 6.11. For this reason, the versions of the Infinity Computer described in the references cited above are similar to the Levi-Civita field implementations, see, e.g., \cite{Shamseddine&Berz(1996)} and \cite{ODE_5}}, and infinitesimal of the same form (\ref{eq:grossnumber}) presented in the decreasing order: $p_0 > p_1 > ... > p_n$. However, according to its description in the papers \cite{ODE_5,Sergeyev_EMS_survey} and patents \cite{Sergeyev_patent}, it can work only under the standard floating-point arithmetic, i.e., only with finite floating-point grossdigits. For this reason, it is difficult to use this methodology to work with neural networks, e.g., using Tensorflow, where the approximations are made using tensors instead of single points $x$.

In this paper, a new technique is proposed to calculate the derivatives of the unknown solution $U(t,x)$ at the initial time $t_0$. The algorithms require using of the software for working with infinitesimal quantities. The required software solution, which can be used not only with floating-point numbers, but with tensors and/or even symbolic expressions, is also presented. The proposed techniques are incorporated with data-driven solution to PDEs using PINNs by Tensorflow. 

The rest of the paper is structured as follows. Section~2 describes the software solution used to deal with infinitesimal numbers and the respective arithmetic. Section~3 presents new differentiation techniques embedded in data-driven solution of PDEs. Section~4 presents the results of the proposed techniques for several real-life problems with known solutions and derivatives. Section~5 presents the results of application of the proposed schemes for data-driven solution by PINNs defined through the DeepXDE software for several real-life problems. Finally, Section~6 concludes the paper. 

\section{Developed software solution for infinitesimal quantities}

There exist several different approaches to work numerically with infinite and infinitesimal numbers: Non-Archimedean fields (e.g., Levi-Civita field, see \cite{Flynn&Shamseddine(2020),Shamseddine&Berz(1996)}), the Infinity Computer (see, e.g., \cite{Sergeyev_EMS_survey}), Dual and Hyper-Dual numbers (see, e.g., \cite{TANAKA201522}), Labeled Sets and Numerosities (see, e.g., \cite{Benci(2003)}), etc. All these methodologies are successfully used in different fields, in particular, in differentiation (see, e.g., \cite{Flynn&Shamseddine(2020),ODE_5,Leonetti&Mukhametzhanov(2022)}). Each of them has its own advantages and limitations, which are briefly described in \cite{Leonetti&Mukhametzhanov(2022)}, and their detailed comparison is out of scope of the present paper. 

All the approaches listed above can be briefly described as follows. A number $C$ representable in a framework for working numerically with infinite and infinitesimal numbers can be defined as 
\begin{equation}
    C = c_0 H_{p_0} + c_1 H_{p_1} + ... c_n H_{p_n},
\end{equation}
where $c_i,~i=0,...,n,$ are standard floating-point numbers, while $H_{p_i},~i=0,...n,$ are finite, infinite, or infinitesimal numbers in the form $H_{p} = h_{1}h_2...h_{k}$, where $k$ is the number of ``axis'' of the approach (e.g., $k = 1$ and $H_p = \go^p$ for the Infinity Computer, $k = 1$ and $H_p = \varepsilon^{p}$ for the Levi-Civita field, $k > 1$ and $H_p = \varepsilon_1^{q_1} \varepsilon_2^{q_2}...\varepsilon_k^{q_k},$ $q_1+q_2+...q_k = p,$ for hyper-dual numbers). For instance, a number $C$ can be representable in the Levi-Civita field as follows:
\begin{equation}
    C = c_0 \varepsilon^{p_0}+c_1\varepsilon^{p_1}+...+c_n\varepsilon^{p_n},
    \label{eq:levi-civita}
\end{equation}
where $c_i\neq 0$ and $p_i,~i=0,1,...,n,$ are finite floating-point numbers, $p_i$ are ordered in increasing order: $p_0 < p_1 < ... < p_n$ (see, e.g., \cite{shamseddine&berz(2010)}). Since in this paper, only derivatives with respect to $t$ are considered, then there is no sense to use hyper-dual numbers, which can be useful for a mixed differentiation due to different ``axis'' of the infinitesimal base. The Infinity Computer with finite powers of $\go$ can also be considered as a particular case of Non-Archimedean fields, e.g., of the Levi-Civita field, since it is easy to see that one can obtain the number in the Infinity Computing form (\ref{eq:grossnumber}) just by fixing $\varepsilon = \go^{-1}$ in (\ref{eq:levi-civita}). However, it should be noted that, in general, the Infinity Computer and the Levi-Civita field are two different approaches (see, e.g., \cite{Sergeyev_EMS_survey}, footnote 10). For this reason, in this paper, the Levi-Civita field is considered as the main used approach for dealing with infinite and infinitesimal numbers and no further considerations refer to the Infinity Computer methodology.

However, since, usually, in data-driven solution, the main data types are tensors instead of single floating-point numbers, then it can be useful to adapt the methodology to work efficiently with tensors. In particular, it is required to represent a tensor $G$ in the following form: 
\begin{equation}
G = G_0 \varepsilon^{p_0} + G_1 \varepsilon^{p_1} + ... + G_n \varepsilon^{p_n},    
\label{eq:tensor_levi-civita}
\end{equation}
where $G_i,~i=0,...,n,$ are tensors of the same dimensions, while $p_i,~i=0,...,n,$ are floating-point numbers written down in a strictly increasing order. Obviously, the tensor $G$ can be also written down in different formats, e.g., as $$G = G_{\varepsilon,0} + G_{\varepsilon,1} + ... + G_{\varepsilon,n},$$ where $G_{\varepsilon,i}$ is a tensor with every element of the form (\ref{eq:levi-civita}). However, in this case, it can be difficult to use efficiently modern frameworks for working with tensors (e.g., Tensorflow \url{https://www.tensorflow.org/} or Pytorch \url{https://pytorch.org/}), since the data types within each tensor are re-defined, while the form (\ref{eq:tensor_levi-civita}) allows to use the advantages of such the frameworks within the Levi-Civita software solution, because it is much simpler to re-define just four main arithmetic operations, than to create another efficient framework for dealing with tensors. 

In this paper, a simplified software solution for working with finite and infinitesimal numbers similar to the Levi-Civita field has been developed. It is based on the Python list comprehensions\footnote{It is clear that faster tools can be used instead of lists, e.g., Numba. Python lists are used for illustration purposes only.} and every quantity $C$ in this framework is represented as follows: 
\begin{equation}
    C = C_0 \varepsilon^0 + C_1 \varepsilon ^1 + ... + C_n\varepsilon^n,
    \label{eq:simplified_levi-civita}
\end{equation}
where $C_i,~i=0,...,n,$ are quantities of the same type and dimensions as $C$ (and can be equal to zero of the same type). Standard finite tensors $C_{finite}$ can be represented in this form simply by $C_0 = C_{finite}$ and $C_i = 0$ for all $i > 0$. In other words, the number $C$ can be represented simply by the vector of the coefficients $[C_0, C_1, ..., C_n]$. Arithmetic operations between the quantities can be simply defined as follows. 

Let $A$ and $B$ be two quantities of the same dimensions in the form (\ref{eq:simplified_levi-civita}) represented by the vectors $[A_0,...,A_n]$ and $[B_0, ..., B_n]$, respectively. Then, $A + B$ is another number of the form (\ref{eq:simplified_levi-civita}) represented by the vector $[A_0 + B_0, A_1 + B_1, ..., A_n + B_n]$. The subtraction is defined in the same way, thus, its definition is omitted. 

The result of multiplication $A\times B$ can also be obtained and represented by the following vector: 
$$[A_0\times B_0, A_0 \times B_1 + A_1\times B_0, A_0\times B_2 + A_1\times B_1 + A_2\times B_0,...,A_n\times B_n],$$
where $A_i \times B_k$ is the element-wise multiplication in the case of tensors. Division can also be defined in a similar way, so its definition is omitted. It should be noted that both multiplication and division should truncate the results after a given number of elements (e.g., after $n-$th element) to preserve the same number of elements within the vector. 

Two important and necessary operations can also be added: multiplication and division by the infinitesimal $\varepsilon^{k},~k>0$. These operations just result in shifting the coefficients to the right (multiplication) or to the left (division). In particular, if $A$ is represented by $[A_0, A_1, ..., A_n]$, then $A\times \varepsilon^{k}$ is represented by $[0,...,0,A_0, A_1,...,A_{n-k}]$, where the first $k$ elements are equal to zero of the same type as $A$ (e.g., zero tensor), and the last elements are truncated to preserve the same dimension of the array. Similarly, $A/\varepsilon^k$ is represented by $[A_k, A_{k+1}, ..., A_{n}, 0, ..., 0]$, where the first $k$ elements have been truncated, while the last $k$ are set to 0. It should be noted that, similarly to fixed-point arithmetic, the presented methodology cannot be used for general purpose problems due to truncations resulting in above mentioned multiplications and divisions. However, this methodology is sufficient for illustration purposes to deal with differentiation problems, where there is no necessity of working with infinite quantities. For more complex problems, there should be constructed a software, e.g., similar to the floating-power simulator from \cite{ODE_5}, adapted to work with tensors. 

It is important to note that there should not arise infinite values in differentiation, if the functions under consideration are differentiable. However, there can be necessary to hold intermediate infinite values: for example, if $g(x) = 1.5x^2,$ then $g(\varepsilon) = 1.5\varepsilon^2$, which is infinitesimal. If $g(x)$ is constructed through complex computational procedures, there can be a situation when $g(x)$ is defined as, e.g., $g(x) = (3x^2/2x^4)/x^{-4}$, where $g(\varepsilon) = (3\varepsilon^{2}/2\varepsilon^4)/\varepsilon^{-4} = 1.5\varepsilon^{-2} / \varepsilon^{-4} = 1.5\varepsilon^2,$ i.e., it is necessary to work with infinite values $\varepsilon^{-4}$ and $\varepsilon^{-2}$ as well. In this case, it is possible to add also several negative powers of $\varepsilon$ in (\ref{eq:simplified_levi-civita}) and to reserve several elements in the respective vector form. In this paper, however, it is supposed that all functions are well-defined and there won't arise situations of this kind just for simplicity. 

\section{Proposed differentiation methods for PDEs}
Let us consider the following nonlinear PDE:
\begin{equation}
    U_t = F(U,U_x,U_{xx},t,x),~U = U(t,x),~t\in [t_0, t_{end}],~x\in D \subset \mathbb{R}^N,
    \label{eq:pde_problem}
\end{equation}
with the following initial conditions:
\begin{equation}
    U(t_0, x) = g(x),~x\in D,
    \label{eq:initial_cond}
\end{equation}
and some boundary conditions. Here, $U(t,x)$ is the unknown solution of the PDE (possibly, multidimensional vector function $U: \mathbb{R}_+\times \mathbb{R}^N \longrightarrow \mathbb{R}^M,$ $M~\geq~1$), $t$ and $x$ are independent variables denoted as time and spatial variables, respectively, $t_0 < t_{end}$ are two real-valued numbers, $D$ is the spatial domain, e.g., a hyperinterval $D = [a,b] \subset \mathbb{R}^N$. Without loss of generality, let us consider $t_0 = 0$ and $t_{end} = T,$ where $T \in \mathbb{R},$ $T > 0$. Let us also consider the univariate case $N = 1$, i.e., the variable $x$ is univariate, just for simplicity and without loss of generality (it will be easy to see that the algorithms are easily generalizable to the case $N > 1$). 

Note also that the equation of the higher than the first order with respect to time can also be transformed into (\ref{eq:pde_problem}) by increasing the dimension and substituting $U_t = V,$ $V_t = F(U,V, U_x, V_x,U_{xx}, V_{xx},t,x)$ (for the second order equation $U_{tt} = F(U,U_x,U_{xx},t,x)$), adding so additional functions and equations (see the Wave equation as an example in the next section). 

In the present paper, the problem of numerical computation of the higher order derivatives $\frac{\partial^k U}{\partial t^k},$ $k = 1,~2,~3,...,$ of $U(t,x)$ at the initial time $t = 0$ is considered. The boundary conditions are not used in the presented algorithms, thus they can be given in any format (e.g., Dirichlet, Neumann, Cauchy, etc.). 

The main idea of the proposed algorithm consists of using the initial solution $U(0,x) = g(x)$ for reconstruction of the Taylor expansion with respect to $t$ of $U(t,x)$ around $t = 0$ simultaneously for all available (fixed) values $x = X$ (e.g., generating $k$ values $x_1, x_2, ..., x_k$ randomly within $D$ and then sending them to $g(x)$ in tensor form $X$ using Tensorflow). The equation (\ref{eq:pde_problem}) at $t = 0$ becomes:
\begin{equation}
    U_t|_{t = 0, x = X} = F(g(X),g'(X),g''(X), 0, X),
    \label{eq:first_derivative}
\end{equation}

If the explicit formulae of $g(x)$ are available, then it is possible to calculate its derivatives with respect to $x$ by any available method (e.g., analytically, symbolically, or by automatic differentiation). Otherwise, if the analytical formulae are not available, the derivatives of $g(x)$ can be obtained by automatic differentiation or by other techniques, e.g., using hyper-dual numbers. 

The equation (\ref{eq:first_derivative}) allows to obtain the first derivative $U_t$ at $t = 0$ and $x = X$. This additional information allows to generate the first-order Taylor expansion of $U(t,X)$ for small $t$. In particular, given the infinitesimal value $h$, the first-order approximation $u_1$ of $U(h,X)$ can be obtained as follows
\begin{equation}
    u_1 = U(0,X) + U_t|_{t=0,x=X}\cdot h = g(X) + F(g(X),g'(X),g''(X),0,X) \cdot h.
    \label{eq:u_1}
\end{equation}

After computation of $u_1$, one can also calculate the value $F(u_1,u_1',u_1'',h,X)$, where $u_1'$ and $u_1''$ denote, respectively, the first and the second derivatives of $u_1$ with respect to $x$ at $t = 0$ and $x = X$. These two derivatives can be again calculated using, e.g., automatic differentiation of $u_1$, since it is calculated from $F(g(X), g'(X), g''(X), 0,X)$, which, in its turn, is calculated from the function $g(X)$, for which the automatic differentiation can be applied. 

Since the value $h$ is represented in the form (\ref{eq:simplified_levi-civita}) (or, more generally, in the form (\ref{eq:levi-civita})), then $u_1$ and $F(u_1,u_1',u_1'',h,X)$ are also written down in the same form. Just for simplicity, hereinafter, the value $h = \varepsilon$ from (\ref{eq:simplified_levi-civita})  represented by the vector $[0,1,0,...,0]$ is used. In this case, $u_1$ is already defined in this form with the tensors $C_0 = g(X)$, $C_1 = F(g(X),g'(X),g''(X),0,X)$ and all $C_i,~i>1,$ are zero tensors, while the value $F(u_1,u_1',u_1'',h,X)$ will be also written down in powers of $\varepsilon$:  
\begin{equation}
    F(u_1,u_1',u_1'',h,X) = F_{1,0}\varepsilon^0 + F_{1,1} \varepsilon^1 + ... + F_{1,n} \varepsilon^n,
    \label{eq:F(u_1)}
\end{equation}
where $F_{1,i}$ are tensors of the same dimension. It is easy to extract the value $F_{1,1}$ from (\ref{eq:F(u_1)}), that is the first-order Lie derivative of $F(...)$ with respect to $t$ at $t = 0$ and $x = X$. Since the function $F$ defines the first-order derivative of $U(t,x)$ at $t = 0, $ $x = X$, then its first-order Lie derivative defines the second-order time derivative of $U$. For this reason, it is possible to re-construct the Taylor expansion of $U$ up to the second derivative obtaining so the value $u_2$ as follows: 
\begin{equation}
    u_2 = U(0,X) + U_t|_{t=0,x=X}\cdot h + \frac{1}{2!} F_{1,1} \cdot h^2 = u_1 + \frac{1}{2!} F_{1,1} \cdot h^2,
    \label{eq:u_2}
\end{equation}
Since $h = \varepsilon$ and represented by $[0,1,0,...,0]$, then $h^2 = \varepsilon^2$ is represented by the vector $[0,0,1,0,...,0]$. Thus, the respective vector representation of $u_2$ can be also obtained: $[g(X), F(g(X), g'(X), g''(X), 0,X), \frac{1}{2} F_{1,1}, 0,...,0]$, i.e., the first two coefficients $C_0$ and $C_1$ are the same as in $u_1$, the coefficient $C_2$ is equal to $\frac{1}{2} F_{1,1}$ and the rest is zero. It should be noted that in (\ref{eq:u_1}) and (\ref{eq:u_2}), the value $h$ is used instead of explicit $\varepsilon$ just to distinguish the infinitesimal \emph{number} $h$ represented by the vector $[0,1,0,...,0]$ in the form (\ref{eq:simplified_levi-civita}) and the infinitesimal base $\varepsilon$, which is used in (\ref{eq:simplified_levi-civita}) and is not explicitly stored in the computer memory. 

Given the value $u_2$, again, the value $F(u_2,u_2',u_2'',h,X)$ can be calculated:
\begin{equation}
    F(u_2,u_2',u_2'',h,X) = F_{2,0}\varepsilon^0 + F_{2,1} \varepsilon^1 + ... + F_{2,n} \varepsilon^n,
    \label{eq:F(u_2)}
\end{equation}
from where, similarly, the coefficient $F_{2,2}$ can also be extracted, which represents the second-order Lie derivative of $F$ divided by $2!$, i.e., the third-order time derivative of $U$ at $t = 0$ and $x = X$. Extracting $F_{2,2}$ from (\ref{eq:F(u_2)}) and multiplying it by $2!$, one can obtain the third-order approximation $u_3$: $u_3 = u_2 + \frac{1}{3!} \cdot 2! \cdot F_{2,2} \cdot h^3 = u_2 + \frac{1}{3} F_{2,2}$. 

It is easy to see that the coefficients $F_{2,0}$ and $F_{2,1}$ coincide with the values $F_{1,0}$ and $f_{1,1}$, respectively, since the first two coefficients in $u_1$ and $u_2$ coincide. Every iteration of the procedure described above improves the order of the approximation of $U(h,X)$ and allows to obtain the required derivatives of $U(t,x)$ at $t = 0$ and $x = X$, so it becomes easy to define the  Algorithm~\ref{alg:diff_alg}. 

\begin{algorithm}[h!]
\caption{Differentiation of $U(t,x)$ from (\ref{eq:pde_problem}) w.r.t. time $t$}
\label{alg:diff_alg}
\begin{algorithmic}[1]
\Statex \textbf{Input} 
\State The problem (\ref{eq:pde_problem}) with the initial condition (\ref{eq:initial_cond}) linked to an automatic differentiation technique (or any other similar) with respect to $x$. 
\State The set (e.g., in tensor form) $X \subset D$ of fixed points within the domain $D$. 
\State The software for working with infinitesimal numbers (e.g., the presented Levi-Civita field solution for the numbers from (\ref{eq:simplified_levi-civita})) and the respective number $n$ of infinitesimal parts to work with from (\ref{eq:simplified_levi-civita}) or from (\ref{eq:levi-civita}) in general. 
\State The infinitesimal $h = 1\varepsilon^1$ in the respective vector form $[0,1,0,...,0]$. 
\Statex \textbf{Result} 
\State The tensor of the derivatives $\frac{\partial^i U}{\partial t^i}~(t,x)~|_{t = 0,x=X},$ $i = 0,1,...,K,$  $K \leq n$. 
\Statex \textbf{Main procedure}
\State $u \gets g(X) + F(g(X), g'(X), g''(X), 0, X)\cdot h$ 
\State $u' \gets g'(X)$, $u'' \gets g''(X)$
\For{$i = 2,3,...,K$}
\State Calculate the first and second order derivatives $u_i'$ and $u_i''$ of the coefficient of $\varepsilon^{i-1}$ at $u$ with respect to $x$ at $x = X$:
\State $u' = u' + u_i' \varepsilon^i$,
\State $u'' = u'' + u_i'' \varepsilon^i$
\State Calculate $F(u, u', u'', h, X)$ up to the coefficient of $\varepsilon^{i-1}$:
\State $F(u_i, u_i', u_i'', h, X) = F_{0}\varepsilon^0 + F_{1}\varepsilon^1 +...+ F_{i-1}\varepsilon^{i-1}$.
\State Save the derivative $\frac{\partial^i U}{\partial t^i}~(t,x)~|_{t = 0,x=X} \gets (i-1)! \cdot F_{i-1}$.
\State Update $u \gets u + \frac{1}{i} \cdot F_{i-1} \cdot h^i$.
\EndFor
\Statex $\bold{Return}$ the obtained derivatives $\frac{\partial^i U}{\partial t^i}~(t,x)~|_{t = 0,x=X}$ or the Taylor expansion $u$ of $U(t,x)$ at $t = h$ and $x = X$.
\end{algorithmic}
\end{algorithm}

\textbf{Remark 1.} Since the values in the form (\ref{eq:simplified_levi-civita}) can be represented by the vectors, then it is not necessary to re-calculate the derivatives of the coefficients with respect to $x$ of the coefficients of $\varepsilon^l$ in $u$ for $l \leq i$: they are already calculated at the previous steps. The vector form of the values $u$, $u'$ and $u''$ is very useful here for this reason (only the last component is differentiated and updated every time). 

\textbf{Remark 2.} There is also produced the Taylor expansion $u = C_0\varepsilon^0 + C_1\varepsilon^1 + ... + C_K\varepsilon^K$ of $U(h,X)$, where $C_0 = g(X)$ and $C_i = \frac{(i-1)!}{i!} F_{i-1}, i = 1,...,K,$ which can also be useful in practice for several reasons. Since, $u$ is represented by the vector, then it is easy to substitute the value $\varepsilon$ by any other finite or infinitesimal value $s$ obtaining so different approximations of $U(s,X)$ for different $s$ close to $0$. For instance, having the Taylor expansion $u$ up to the 10-th derivative, it is possible to calculate the value $U(0.1,X)$ just multiplying element-wise the vector $[C_0, C_1,...,C_K]$ by the vector $[1, 0.1, 0.1^2,...,0.1^K]$. This issue allows to obtain high precision approximations of $U$ for $t > 0$ generating so additional data values and improving the performance of the respective PINNs. In the next sections, there will be presented several examples of such the improvements. 

\textbf{Remark 3.} The accuracy can be lower for the highest derivatives (e.g., of the order 10), if the absolute value of the derivative is high (see, e.g., the Wave equation in the next section with large values of the parameter $C$). This occurs due to the limitations of the floating-point arithmetic and due to the use of factorials in the computation of the derivatives: the $i-$th derivative of $U$ is obtained multiplying the coefficient $F_{i-1}$ by $(i-1)!$ (see the lines 11-12 of the Algorithm~\ref{alg:diff_alg}). This issue will be studied in the future work in order to develop a more efficient algorithm of computations without involving the factorials and, thus, without loosing the accuracy even for the highest derivatives. 

It should be noted that the proposed techniques are referenced with the initial condition $U(0,x) = g(x)$. However, in data-driven solution, the initial condition function $g(x)$ can unknown. The set of values $U(0,x_i) = g_i$ can be given instead. In this case, the algorithms presented in this paper can also be used, e.g., as follows. First, the function $g(x)$ can be approximated by another differentiable neural network using both initial and boundary conditions. Then, the derivatives with respect to $x$ (required to calculate the derivatives with respect to $t$) can be calculated, e.g., using Tensorflow Autodiff, allowing so to calculate the time derivatives as well (time derivatives cannot be calculated in the same way by Autodiff, since the initial conditions are defined for only one time value $t = 0$). 

\section{Numerical experiments on benchmarks with known solutions}
Numerical experiments have been carried out in Jupyter notebooks using Google Colab Pro. The software has been written in Python 3.7.14. Autodiff of Tensorflow (version 2.8.2) has been used for the computation of the derivatives of $U$ with respect to $x$. 

For each test problem, the derivatives up to the order $10$ have been calculated $\frac{\partial^i U}{\partial t^i}(t,x),$ $i = 0,1,...,10,$ at the time $t_0 = 0$ on $50$ different points $x_i$ chosen uniformly within the interval $D$ excluding the points, where the solution $U(t_0,x)$ is too small, just to avoid divisions by 0 in the relative error computations. All computations have been performed in tensor form for all 50 points simultaneously, i.e., $X = \mbox{tensorflow.Variable}(x, dtype=\mbox{'float64'})$, where $x$ is the array of 50 different points. 

Once all 10 derivatives have been obtained, the Normalized Root Mean Square Error (NRMSE) has been calculated for each derivative with respect to the exact derivatives (known analytically for each test problem): 
\begin{equation}
    NRMSE(y_{true}, y_{approx}) = \frac{1}{N} \frac{||y_{true} - y_{approx}||_2}{y_{true}^{MAX} - y_{true}^{MIN} + 1},
\end{equation}

As it has been already mentioned before, the error of the computation of the highest derivatives cannot be a unique metric to determine the quality of the obtained information. For this reason, the Taylor expansion at $t =  0.01,~0.05,~0.1,$  has been also calculated using the obtained value $u$ being the Taylor approximation of $U(\varepsilon,X)$ around $t = 0$ and $x = X$. The NRMSE has been also calculated with respect to the exact value of the solution (the exact solution $U(t,x)$ is also known for each test problem).

\subsection{Heat equation}
Let us consider the following PDE from \cite{Raissi:et:al.(2018)}: 
\begin{equation}
    U_t = \alpha U_{xx},~x\in [0,L], t\in [0,1],
    \label{eq:heat}
\end{equation}
with the periodic initial condition: 
\begin{equation}
    U(0,x) = \sin{(\frac{n\pi x}{L})},~x\in [0,L],~n=1,2,...,
\end{equation}
where $\alpha = 0.4 $ is the thermal diffusivity constant, $L = 1$ is the length of the bar, $n = 1$ is the frequency of the sinusoidal initial conditions. The Dirichlet boundary conditions are given as follows: 
\begin{equation}
    U(t,0) = U(t,L) = 0,~t \geq 0.
\end{equation}
The exact solution of (\ref{eq:heat}) is known: 
\begin{equation}
    U(t,x) = e^{\frac{-n^2\pi^2\alpha t}{L^2}} \cdot \sin{(\frac{n\pi x}{L})},
\end{equation}
from where it is also possible to calculate the exact derivatives of $U$ with respect to time:
\begin{equation}
    \frac{\partial^i U}{\partial t^i}(t,x) = (\frac{-n^2\pi^2\alpha}{L^2})^i \cdot U(t,x).
\end{equation}

The Heat equation (\ref{eq:heat}) is quit simple, but allows to study well the proposed algorithm. The respective values of the NRMSE for the first 10 derivatives are presented in Table~\ref{tab:benchmarks_rmse} (column Heat). As it can be seen, the obtained errors are of the same order as the machine precision, giving so a lot of additional information about the dynamics of the model. 

\begin{table}[h!]
    \centering
    \begin{tabular}{|c|c|c|c|c|}
    \hline
         \multirow{2}{*}{Diff. order}& \multirow{2}{*}{Heat} & Diffusion & Diffusion & \multirow{2}{*}{Wave}\\
         &&Derivatives&Taylor coeff.&\\
         \hline
       0 & 0 & 0& 9.23e-30  &0\\
1 & 1.70e-16 & 2.33e-16 & 2.33e-16 & 0\\
2 & 3.00e-16 & 2.89e-15 & 1.45e-15 &2.37e-16\\
3 & 3.55e-16 & 4.69e-14 & 7.82e-15 &0\\
4 & 3.98e-16 & 4.37e-13 & 1.82e-14 & 3.10e-16\\
5 & 4.23e-16 & 3.97e-12 & 3.31e-14 &0\\
6 & 9.77e-16 & 2.46e-11 & 3.42e-14 & 3.32e-16\\
7 & 6.88e-16 & 4.20e-10 & 8.34e-14 &0\\
8 & 6.47e-16 & 2.65e-09 & 6.57e-148 & 4.01e-16\\
9 & 8.05e-16 & 3.06e-08 & 8.42e-14 &0\\
10 & 6.61e-16 & 3.61e-07 & 9.94e-14 & 4.31e-16\\
\hline
    \end{tabular}
    \caption{The obtained NRMSE values for the first 10 derivatives of the Heat equation, Diffusion Equation (including the Taylor coefficients NRMSE) and Wave Equation. All errors are of the same order than the machine precision ($\approx 2.2\cdot 10^{-16}$), except the Diffusion equation, for which the Taylor coefficients are accurate instead.}
    \label{tab:benchmarks_rmse}
\end{table}

\begin{figure}[h!]
    \centering
    \begin{minipage}{0.32\textwidth}
    \includegraphics[scale=0.08]{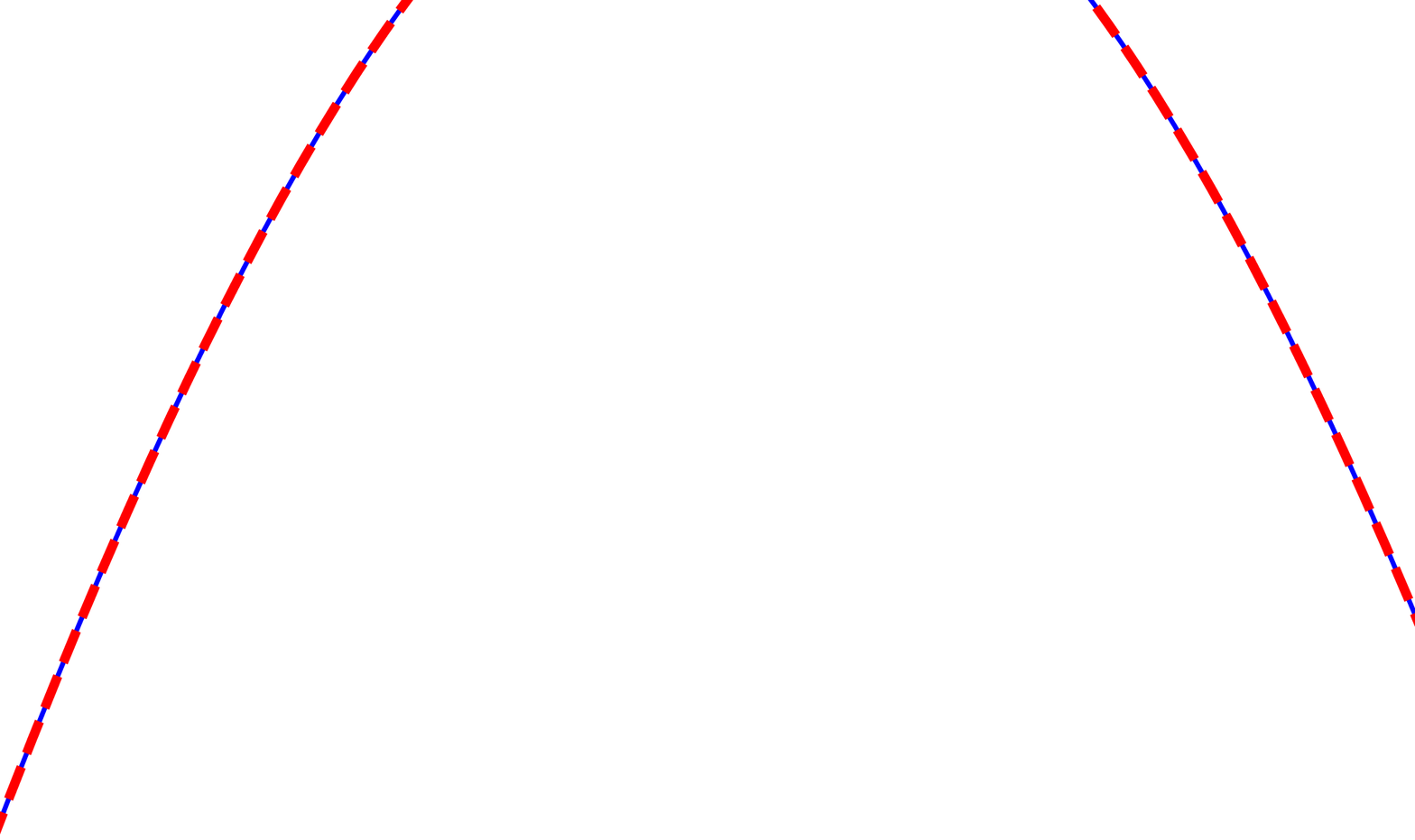} \\ \centering (a)
    \end{minipage}
    \begin{minipage}{0.32\textwidth}
    \includegraphics[scale=0.08]{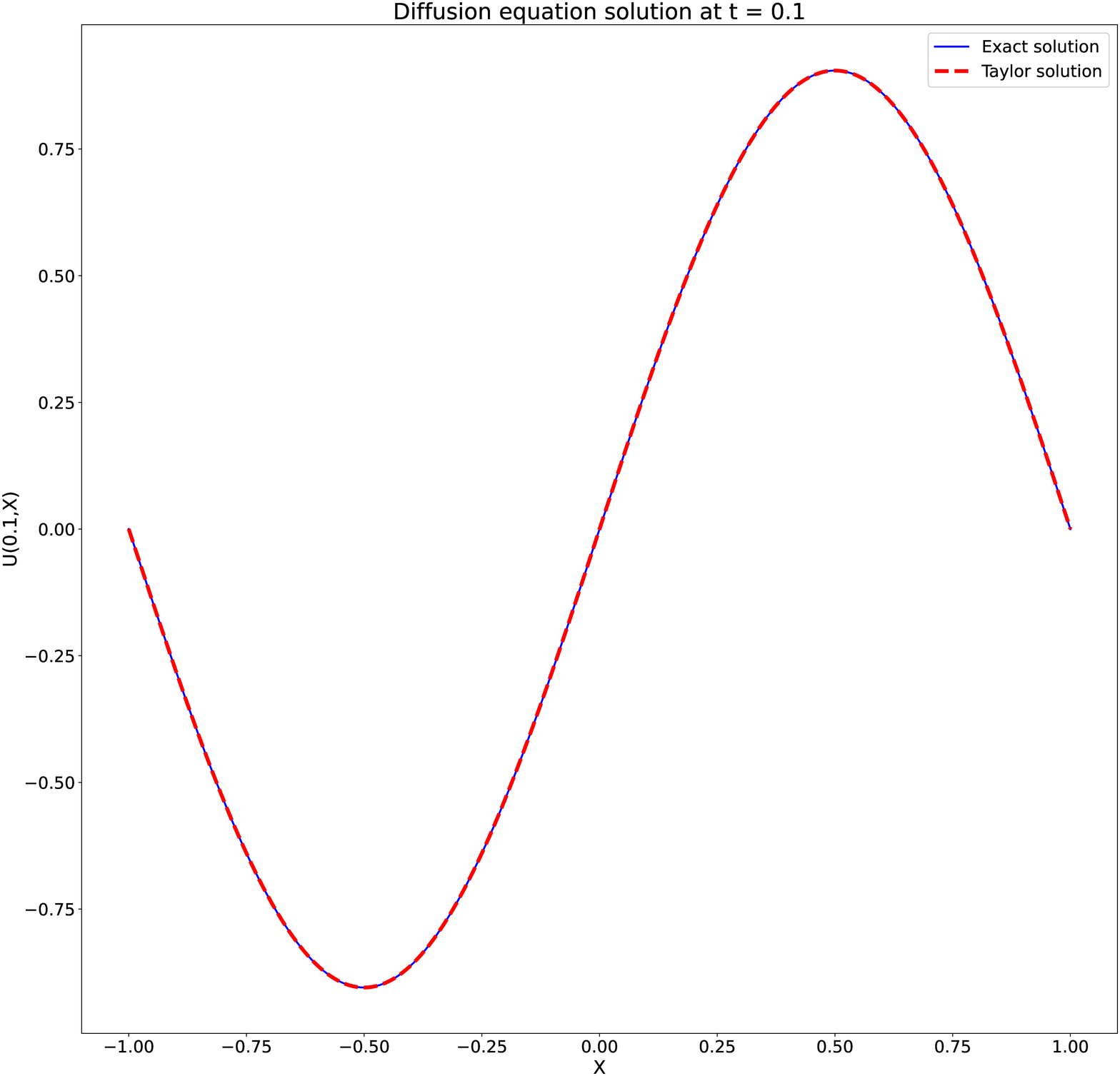} \\\centering (b)
    \end{minipage}
    \begin{minipage}{0.32\textwidth}
    \includegraphics[scale=0.08]{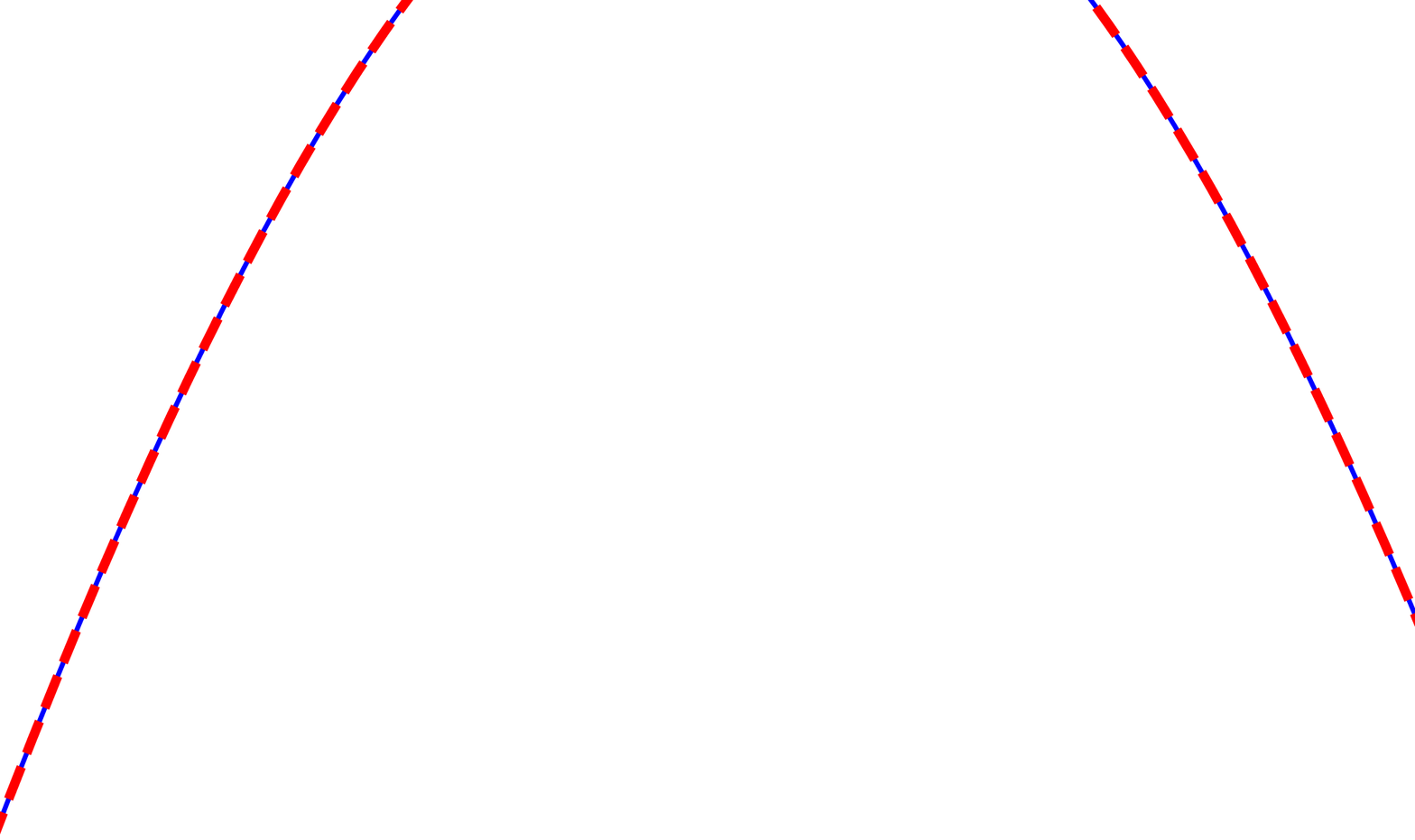} \\ \centering(c)
    \end{minipage}
    \caption{The obtained Taylor approximation using 10 derivatives (red) and the exact solution (blue) of the Heat, Diffusion and Wave equations at $t = 0.1$. Both the solutions coincide up to a very small error for all three test problems.}
    \label{fig:benchmarks}
\end{figure}

The NRMSE values of the Taylor expansion at the time $t_1 = 0.01$, $0.05$, $0.1$ and the same $X$ have been also calculated with respect to the exact solution. One can see from Table~\ref{tab:benchmarks_rmse_tay} (column Heat) that the results have been calculated with a high accuracy (e.g., the NRMSE at $t_1 = 0.01$ is smaller than the machine precision, while the one at $t_1 = 0.05$ is of the same order than the machine precision) giving so additional information about the solution $U(t,x)$, which was unknown a priori. The graph of both the obtained Taylor solution and the exact analytical one at $t = 0.1$ for $500 $ uniform points $x\in D$ is presented in Figure~\ref{fig:benchmarks}(a).

\begin{table}[h!]
    \centering
    \begin{tabular}{|c|c|c|c|}
    \hline
         $t_1$ & Heat & Diffusion & Wave\\
         \hline
        0.01 & 7.75e-17& 3.03e-17 & 1.08e-16\\
        0.05 & 3.63e-16& 3.25e-17 & 1.01e-16\\
        0.1 & 6.63e-13 & 4.22e-17& 2.16e-15\\
\hline
    \end{tabular}
    \caption{The obtained NRMSE values for the Taylor expansion of the Heat, Diffusion and Wave equations at $t_1 = 0.01,~0.05,~0.1$. The obtained approximations have a high accuracy being so useful candidates for an additional data for PINNs.}
    \label{tab:benchmarks_rmse_tay}
\end{table}

\subsection{Diffusion equation}
Let us consider the following PDE from \cite{Raissi:et:al.(2018)}: 
\begin{equation}
    U_t = U_{xx} - e^{-t}(\sin{(\pi x)} - \pi^2 \sin{(\pi x)}),~x\in [-1,1], t\in [0,1],
    \label{eq:diffusion}
\end{equation}
with the initial condition: 
\begin{equation}
    U(0,x) = \sin{(\pi x)},~x\in [-1,1].
\end{equation}
The Dirichlet boundary conditions are given as follows: 
\begin{equation}
    U(t,-1) = U(t,1) = 0,~t \geq 0.
\end{equation}
The exact solution of (\ref{eq:heat}) is known: 
\begin{equation}
    U(t,x) = e^{-t} \cdot \sin{(\pi x)},
\end{equation}
from where it is also possible to calculate the exact derivatives of $U$ with respect to time:
\begin{equation}
    \frac{\partial^i U}{\partial t^i}(t,x) = (-1)^i \cdot U(t,x).
\end{equation}

The equation (\ref{eq:diffusion}) is a bit more complicated, than the previous one, due to the presence of the arguments $t$ and $x$ explicitly in the function $F$. For this reason, it also becomes an important illustrative example for the proposed algorithm. The respective values of the NRMSE for the first 10 derivatives are presented in Table~\ref{tab:benchmarks_rmse} (columns Diffusion Derivatives and Diffusion Taylor coeff.). As it can be seen, the accuracy decreases with the order of the derivatives. This can be easily explained by the inefficient use of the factorials and by the limitations of the floating-point arithmetic. 

The latter fact can be substantiated by the following computations. Let us study directly the coefficients $C_i$ of the Taylor expansion $u$ from the Algorithm~\ref{alg:diff_alg} without multiplying them by the respective factorials, i.e., for $h = \varepsilon$ in Algorithm~\ref{alg:diff_alg}: $u = C_0\varepsilon^0 + C_1\varepsilon^1 + ... + C_n\varepsilon^n$, where $C_i$ is the tensor of the respective dimension (equal to the $i-$th time derivative divided by $i!$ at $t = 0$ and $x = X$).  Let us then calculate the NRMSE for each coefficient $C_i$ with respect to the exact values $\frac{1}{i!}\frac{\partial^i U}{\partial t^i}(0,X)$ (all these computations are made using decimal type to maintain the precision). The results are presented in the last column of Table~\ref{tab:benchmarks_rmse}.

One can see that the accuracy in computation of the derivatives was lost just due to the limitations of the floating-point numbers, while the accuracy in direct coefficients of the Taylor expansion is much higher. The latter fact can also be substantiated by Table~\ref{tab:benchmarks_rmse_tay} (column Diffusion), where the NRMSE values of the Taylor expansion at the time $t_1 = 0.01$, $0.05$, $0.1$ and the same $X$ have been also calculated  with respect to the exact solution. One can see from Table~\ref{tab:benchmarks_rmse_tay} that the results have been calculated up to machine precision at all three points giving so a lot of additional information about the solution $U(t,x)$, which was unknown a priori also in this case. The graph of both the obtained Taylor solution and the exact analytical one at $t = 0.1$ for $500 $ uniform points $x\in D$ is presented in Figure~\ref{fig:benchmarks}(b).

\subsection{Wave propagation equation}
Let us consider the following second-order PDE from \cite{WANG2021113938}: 
\begin{equation}
    U_{tt} = C^2 \cdot U_{xx},~x\in [0,1], t\in [0,1],~C > 0,
    \label{eq:wave}
\end{equation}
with the initial condition: 
\begin{equation}
\begin{matrix}
    U(0,x) & = & \sin{(\pi x)} + \sin{(A\pi x)},~x\in [0,1],~A > 0,\\
    U_t(0,x) &= &0.
\end{matrix}
\end{equation}
The Dirichlet boundary conditions are given as follows: 
\begin{equation}
    U(t,0) = U(t,1) = 0,~t \geq 0.
\end{equation}
The exact solution of (\ref{eq:heat}) is known: 
\begin{equation}
    U(t,x) = \sin{(\pi x)\cos{(C\pi t)}} + \sin{(A\pi x)}\cos{(AC\pi t)},
\end{equation}
from where it is also possible to calculate the exact derivatives of $U$ with respect to time:
\begin{equation}
    \frac{\partial^i U}{\partial t^i}(t,x) =
    (C\pi)^i\sin{(\pi x)}f_1(t) + (AC\pi)^i\sin{(A\pi x)} f_2(t),
\end{equation}
where 
\begin{equation}
\begin{cases}
     f_1(t) = \cos{(C\pi t)},~f_2(t)= \cos{(AC\pi t)}, & \mbox{if } i\%4 = 0,\\
    f_1(t) = - \sin{(C\pi t)},~f_2(t) = -\sin{(AC\pi t)}, & \mbox{if } i\%4 = 1,\\
    f_1(t) = - \cos{(C\pi t)},~f_2(t) = - \cos{(AC\pi t)}, & \mbox{if } i\%4 = 2,\\
    f_1(t) = \sin{(C\pi t)},~f_2(t) = \sin{(AC\pi t)}, & \mbox{if } i\%4 = 3,
    \end{cases}
\end{equation}

The main difference of this problem with respect to the previous ones consists of the order of the equation (\ref{eq:wave}). Since the Algorithm~\ref{alg:diff_alg} has been proposed for the first order equations only, then it should be transformed increasing so the dimension of the problem $Y_t = F(Y,Y_x, Y_{xx}, t,x),$ where $Y = (U,V),$ $V = U_t$, i.e. $Y: \mathbb{R}_+\times \mathbb{R}^N \longrightarrow \mathbb{R}^2$:
\begin{equation}
\begin{cases}
     U_t = V,\\
     V_t = C^2\cdot U_{xx},
    \end{cases}
\end{equation}
and the respective initial and boundary conditions. 

Here, the values of the parameters $A$ and $C$ are fixed equal both to $1$ just for simplicity. The obtained results are presented in Table~\ref{tab:benchmarks_rmse} (column Wave). One can see that the derivatives have been calculated with a high accuracy (of the same order than the machine precision). It should be noted that the errors for the odd orders are equal to zero, since $f_1(t)$ and $f_2(t)$ in (\ref{eq:wave}) are both equal to zero at $t = 0$.

The NRMSE values of the Taylor expansion at the time $t_1 = 0.01$, $0.05$, $0.1$ and the same $X$ have been also calculated with respect to the exact solution. One can see from Table~\ref{tab:benchmarks_rmse_tay} that the results have been again calculated with a high accuracy giving so a lot of additional information about the unknown solution $U(t,x)$. The graph of both the obtained Taylor solution and the exact analytical one at $t = 0.1$ for $500 $ uniform points $x\in D$ is presented in Figure~\ref{fig:benchmarks}(c).

\section{Applications for PINNs}

As it can be seen from the previous section, the proposed algorithm allows to obtain a new information about the dynamics of the model described by the PDE. Let us study how this new information allows to improve the quality of the PINNs for data-driven solution of the PDEs. The schemes proposed in this section are not unique, but just the simplest ones allowing already to increase the accuracy of the solution. 

First of all, the PINNS from the well-known DeepXDE library written in Python with Tensorflow (see, \url{https://github.com/lululxvi/deepxde} and \cite{deepxde}) have been used to solve the PDEs. The library DeepXDE has been chosen for its simplicity, since it allows to define different types of the initial and boundary conditions without changing the code of the PINNs. However, the proposed techniques can also be used for other libraries, e.g., from \cite{nature_pinns}. 

For each test problem, the first $N$ derivatives at $t = 0$ for 100 randomly chosen points $x$ have been calculated. Then, the Taylor approximations using the obtained derivatives have been calculated at $t = 0.01,$ $0.02,$ $0.03,$ $0.04$ and $0.05$ for the same values $x$ obtaining so $500$ additional approximations of $U(t,x)$. The obtained additional values have been added as the Point set conditions into the DeepXDE solver. Comparison with two versions with and without the additional conditions has been then performed. 

For each test problem, forward PINNs have been generated with the following parameters: $5$ inner layers with $64$ neurons each; $8000$ collocation (domain) points for minimizing the loss function on the PDE, $400$ points on the boundary conditions, $800$ points at $t = 0$ by initial conditions. First, each model has been compiled using Adam optimizer with $15000$ iterations (epochs). Then, the training continued with L-BFGS algorithm until stopping criterion. The $L_2$ norm has been calculated on $50$ random points $x \in D$ at some time $t_1$ (chosen similarly for all test problems). All remaining parameters and options have been set to their default values. Loss functions have been set directly using the internal DeepXDE options (deepxde.icbc classes for the initial and boundary conditions) without setting the formulae explicitly.  

For each test problem, the PINN with the above mentioned parameters has been compiled and trained. 
Then, in order to study the usefulness of the proposed techniques, the PINN with the same parameters and additional conditions has been also defined. In particular, $N$ first derivatives $\frac{\partial^i U}{\partial t^i} (0,X)$ have been calculated at $100$ randomly chosen points $X$. Then, the Taylor expansions at $t = 0.01,$ $0.02,...,$ $0.05,$ have been constructed using the calculated derivatives. The obtained $500$ Taylor values have been added as additional conditions to the model as deepxde.icbc.PointSetBC. After that, the model has been compiled and trained in the same way with the same parameters. The results are presented for each test problems in the following subsections. The Python code for the main PINN implementation and visualizations has been taken from \url{https://deepxde.readthedocs.io/en/latest/demos/pinn_forward.html}. The obtained results and models configurations can be also found in \url{https://github.com/maratsmuk/DiffPDE_PINNs} for transparency and illustration aims only.

\subsection{Burgers Equation}
The first problem under consideration was the well-known Burgers' equation (see, e.g., \cite{BASDEVANT198623}). 
\begin{equation}
    \frac{\partial U}{\partial t} = - U \frac{\partial U}{\partial x} + \nu \frac{\partial^2 U}{\partial x^2},
    \label{eq:burgers}
\end{equation}
where $\nu = \frac{0.01}{\pi}$ is the constant and with the following initial and boundary conditions:
\begin{equation}
\begin{matrix}
    U(t,-1) = U(t,1) = 0,\\
    U(0,x) = -\sin (\pi x),\\
    x \in [-1,1],~t\in[0,1],~U = U(t,x) \in \mathbb{R}.
    \end{matrix}
    \label{eq:burgers_icbc}
\end{equation}

\begin{figure}
    \centering
    \hspace{-0.12\textwidth}
    \begin{minipage}{0.45\textwidth}
    \includegraphics[scale=0.15]{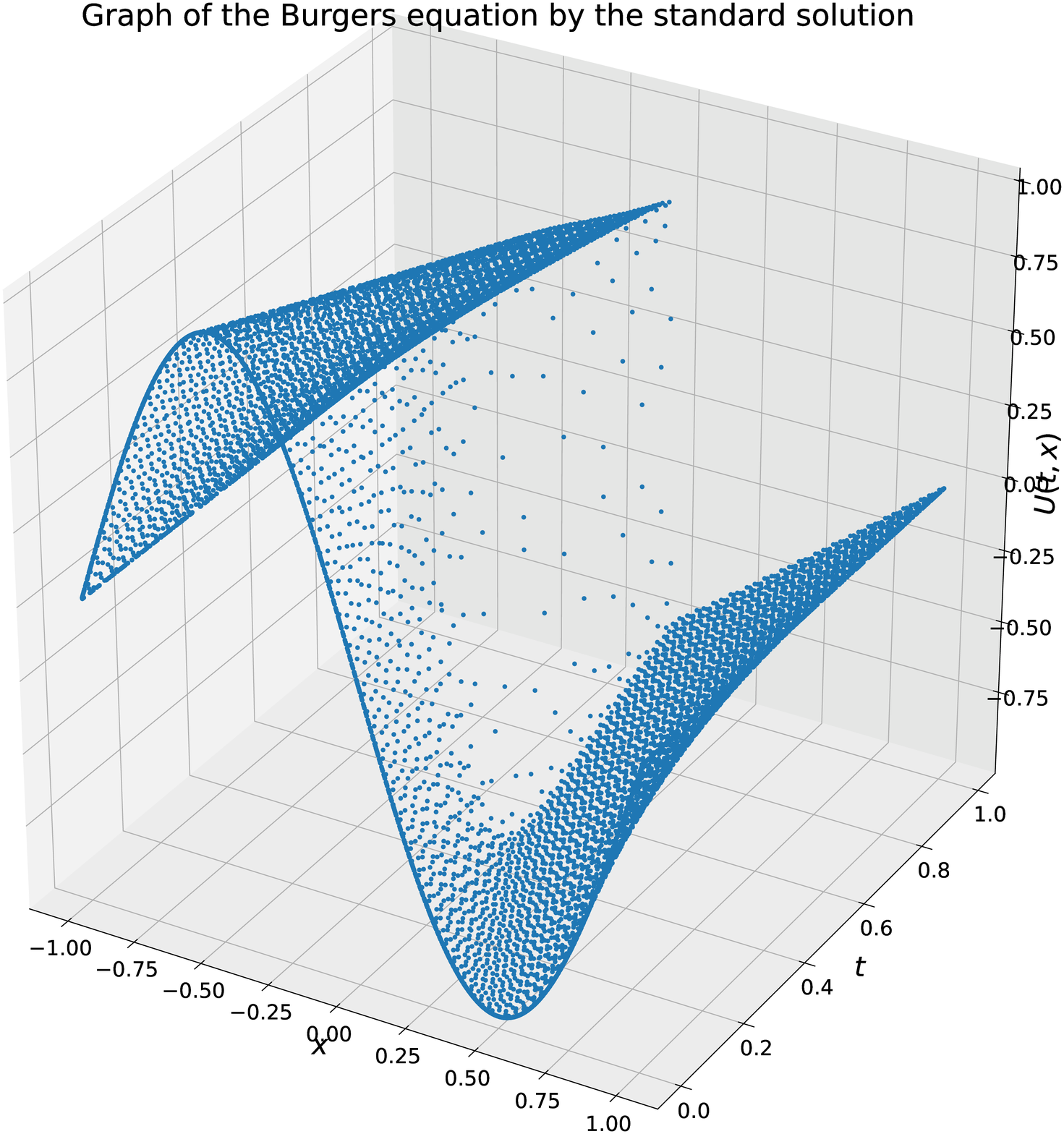}\\
    \begin{center}
        \vspace*{-0.25\textwidth} \hspace*{0.35\textwidth}(a)
    \end{center}
    \end{minipage}
    \hfill
    \begin{minipage}{0.45\textwidth}
    \hspace{-0.2\textwidth}
    \includegraphics[scale=0.15]{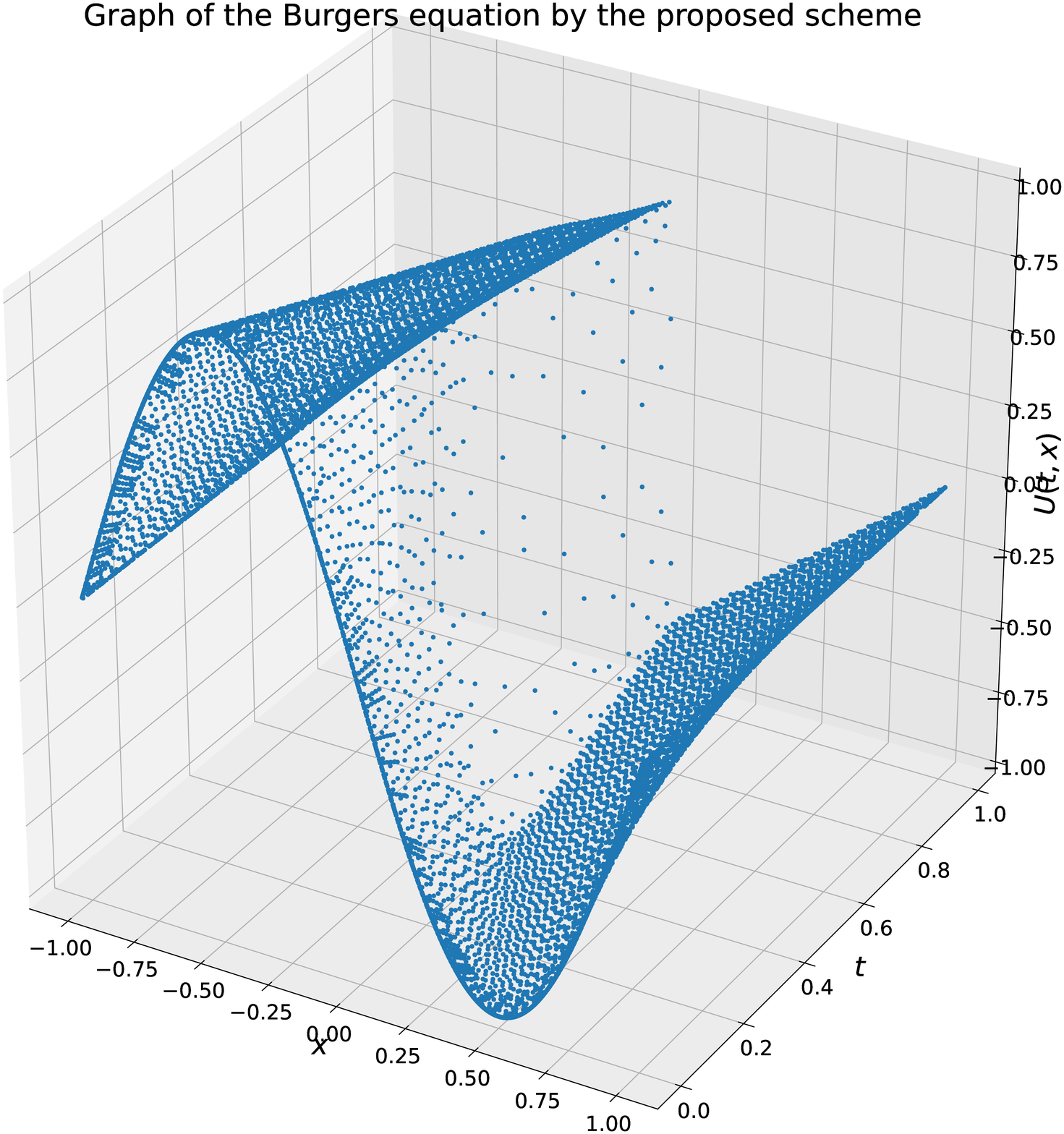} \\
    \begin{center}
        \vspace*{-0.25\textwidth} \hspace*{0.1\textwidth}(b)
    \end{center}
    \end{minipage}
    \caption{The obtained solutions by PINNs using the standard model (left) and with additional conditions by the proposed scheme (right) of the Burgers' equation. }
    \label{fig:burgers}
\end{figure}

First, the Dirichlet boundary conditions and standard initial conditions have been generated. The best train loss of the standard model has been obtained: $8.13e-08$.
The $L_2$ norm of the obtained solution at $t = 0.1$ and $50$ random points $x$ has also been calculated: $4.183e-05$ with the mean residual $0.0002$. The plot of the obtained solution is presented in Figure~\ref{fig:burgers}(a). 

Then, the first $7$ derivatives have been calculated at $t = 0$. The obtained Taylor approximations have been calculated for $100$ random points $x$ at $t = 0.01,$ $0.02,...,$ $0.05$. The obtained $500$ values have been added as additional conditions to the model, which was re-initialized, compiled and trained with the same parameters. The final obtained loss has become $8.38e-08$. The obtained $L_2$ norm at $t = 0.1$ has become better: $3.612e-05$ with the mean residual $0.000227$. The plot of the obtained solution is also presented in Figure~\ref{fig:burgers}(b). The obtained $L_2$ norm has become better with almost the same number of epochs (31321 epochs against 31347 for the standard scheme) showing so a promising behavior of the proposed scheme. 

\subsection{Allen-Cahn Equation}
The second problem under consideration was the well-known Allen-Cahn equation (see, e.g., \cite{RAISSI2019686}):
\begin{equation}
\begin{matrix}
    \frac{\partial U}{\partial t} =  d \frac{\partial^2 U}{\partial x^2} + 5 (U - U^3),
    \end{matrix}
    \label{eq:allen-cahn}
\end{equation}
where $d = 0.0001$ and with the following initial and boundary conditions:
\begin{equation}
\begin{matrix}
    U(t,-1) = U(t,1),\\
    U_x(t,-1) = U_x(t,1),\\
    U(0,x) = x^2 \cos (\pi x),\\
    x \in [-1,1],~t\in[0,1],~U = U(t,x) \in \mathbb{R},
    \end{matrix}
    \label{eq:allen-cahn_icbc}
\end{equation}

\begin{figure}
    \centering
    \hspace{-0.12\textwidth}
    \begin{minipage}{0.45\textwidth}
    \includegraphics[scale=0.15]{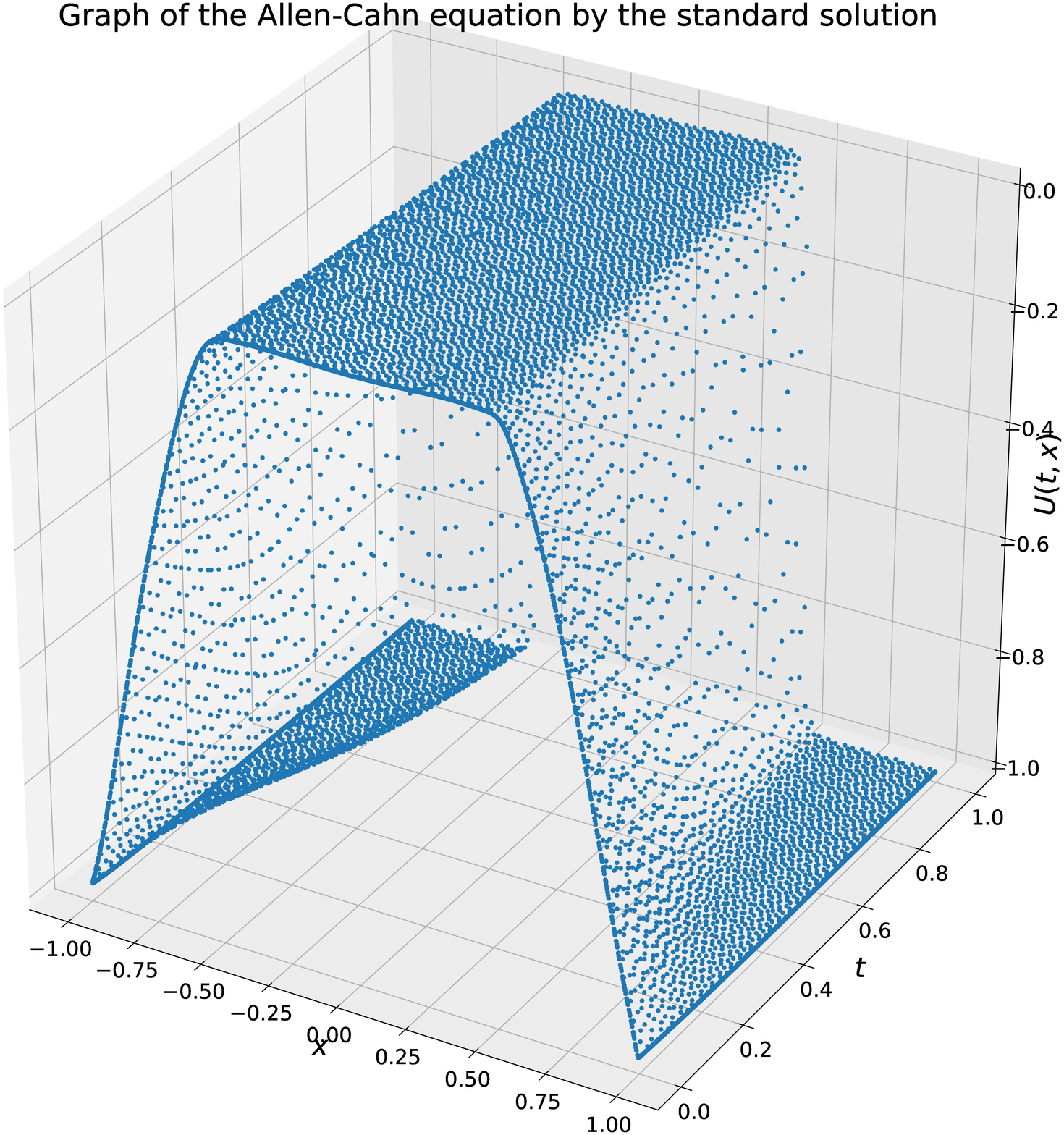}\\
    \begin{center}
        \vspace*{-0.25\textwidth} \hspace*{0.35\textwidth}(a)
    \end{center}
    \end{minipage}
    \hfill
    \begin{minipage}{0.45\textwidth}
    \hspace{-0.2\textwidth}
    \includegraphics[scale=0.15]{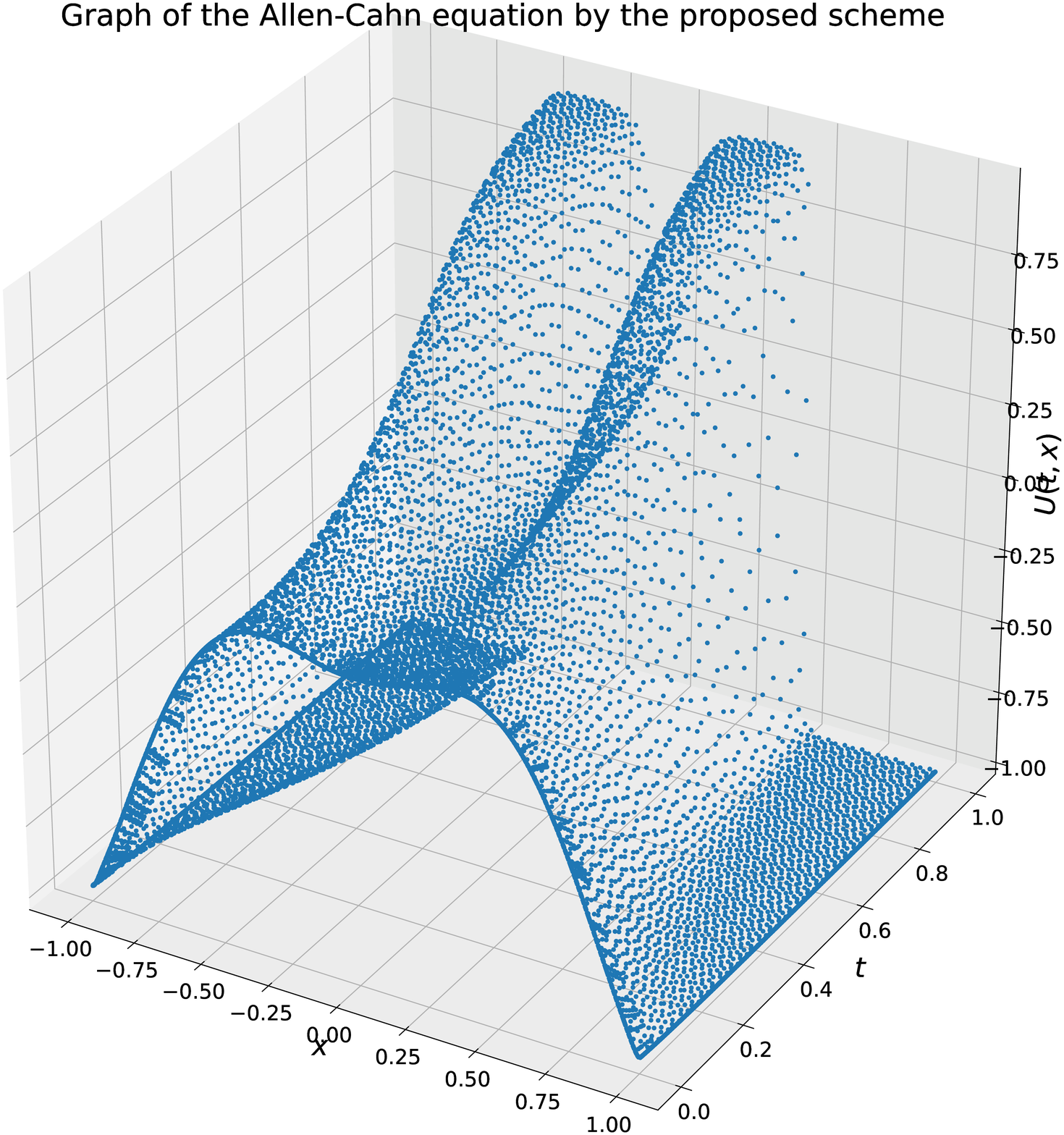} \\
    \begin{center}
        \vspace*{-0.25\textwidth} \hspace*{0.1\textwidth}(b)
    \end{center}
    \end{minipage}
    \caption{The obtained solutions by PINNs using the standard model (left) and with additional conditions by the proposed scheme (right) of the Allen-Cahn equation. One can see that the standard model failed  to approximate the solution correctly, while additional data by the proposed scheme allowed to improve the accuracy a lot.}
    \label{fig:allen-cahn}
\end{figure}

First, the Periodic boundary conditions for both $U$ and $U_x$ and standard initial conditions have been generated. The best train loss of the standard model has been obtained: $6.45e-04$.
The $L_2$ norm of the obtained solution at $t = 0.1$ and $50$ random points $x$ has also been calculated: $0.0584$ with the mean residual $0.0073$. The plot of the obtained solution is presented in Figure~\ref{fig:allen-cahn}(a). 

Then, the first $7$ derivatives have been calculated at $t = 0$. The obtained Taylor approximations have been calculated for $100$ random points $x$ at $t = 0.01,$ $0.02,...,$ $0.05$. The obtained $500$ values have been added as additional conditions to the model, which was re-initialized, compiled and trained with the same parameters. The final obtained loss has become $7.37e-06$. The obtained $L_2$ norm at $t = 0.1$ has become much better: $0.0017$ with the mean residual $0.00146$. The plot of the obtained solution is also presented in Figure~\ref{fig:allen-cahn}(b). One can see from Figure~\ref{fig:allen-cahn} that the standard model was not able to approximate the solution $U(t,x)$ correctly within the given iterations. Moreover, the overfitting has clearly arisen around the initial solution. The proposed scheme, instead, has allowed to obtain a much better solution without overfitting, due to additional data obtained from the calculated derivatives. 

\subsection{Nonlinear Schrodinger Equation}
Finally, the well-known Nonlinear Schrodinger equation has also been studied (see, e.g., \cite{RAISSI2019686}): 
\begin{equation}
    i \frac{\partial h}{\partial t} =  -0.5 \frac{\partial^2 h}{\partial x^2} - |h|^2 h,
    \label{eq:schrodinger}
\end{equation}
where  is a complex valued function and with the following initial and boundary conditions:
\begin{equation}
\begin{matrix}
    h(0,x) = 2 sech(x),\\
    h(t,-5) = h(t,5),\\
    h_x(t,-5) = h_x(t,5),\\
    x \in [-5,5],~t\in[0,\pi/2],~h(t,x) = U(t,x) + iV(t,x) \in \mathbb{C}.
    \end{matrix}
    \label{eq:schrodinger_icbc}
\end{equation}

In order to train the model, the function $h$ has been considered as a two-dimensional function $h = (U,V)$, which was represented by the neural network with two outputs. The equation (\ref{eq:schrodinger}) has been defined as a system of two PDEs for the real and imaginary parts, respectively, allowing one to calculate the derivatives $\frac{\partial^i h}{\partial t^i} (0,X)$ and to compile and train the respective models. 

\begin{figure}
    \centering
    \hspace{-0.12\textwidth}
    \begin{minipage}{0.45\textwidth}
    \includegraphics[scale=0.15]{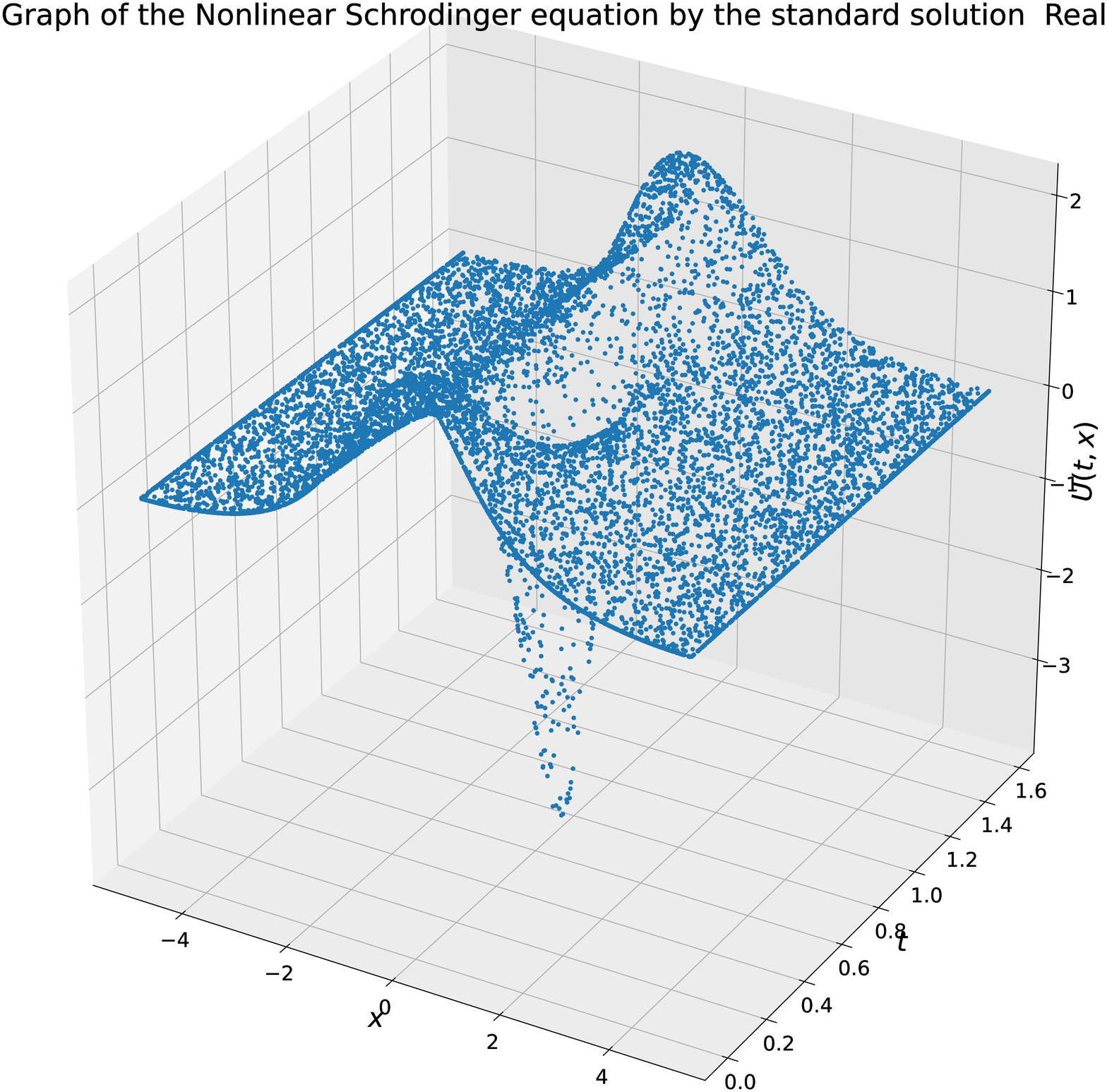}\\
    \begin{center}
        \vspace*{-0.25\textwidth} \hspace*{0.35\textwidth}(a)
    \end{center}
    \end{minipage}
    \hfill
    \begin{minipage}{0.45\textwidth}
    \hspace{-0.2\textwidth}
    \includegraphics[scale=0.15]{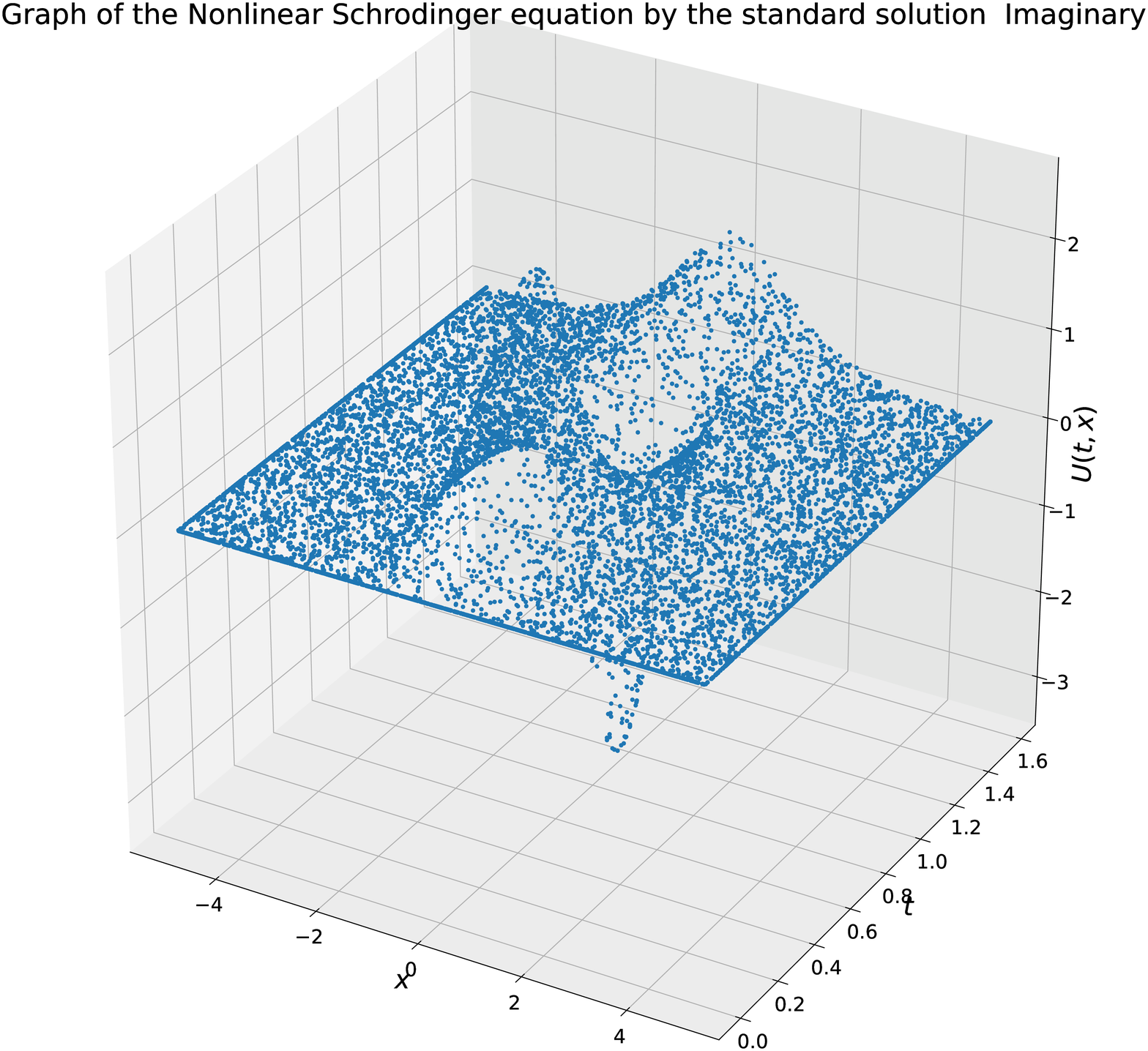} \\
    \begin{center}
        \vspace*{-0.25\textwidth} \hspace*{0.1\textwidth}(b)
    \end{center}
    \end{minipage}
    \vfill
    \centering
    \hspace{-0.12\textwidth}
    \begin{minipage}{0.45\textwidth}
    \includegraphics[scale=0.15]{NonlinearSchrodinger0False_sol.eps}\\
    \begin{center}
        \vspace*{-0.25\textwidth} \hspace*{0.35\textwidth}(c)
    \end{center}
    \end{minipage}
    \hfill
    \begin{minipage}{0.45\textwidth}
    \hspace{-0.2\textwidth}
    \includegraphics[scale=0.15]{NonlinearSchrodinger1False_sol.eps} \\
    \begin{center}
        \vspace*{-0.25\textwidth} \hspace*{0.1\textwidth}(d)
    \end{center}
    \end{minipage}
    \caption{The obtained solutions by PINNs using the standard model (top, figures (a) and (b)) and with additional conditions by the proposed scheme (bottom, figures (a) and (b)) of the nonlinear Schrodinger equation.}
    \label{fig:schrodinger}
\end{figure}

First, the Periodic boundary conditions for both $h$ and $h_x$ and standard initial conditions have been generated separately for the real and imaginary parts of $h$. The best train loss of the standard model has been obtained: $3.05e-06$.
The $L_2$ norm of the obtained solution at $t = 0.1021$ (the closest to $0.1$ time value available in the data set) and $50$ random points $x$ has also been calculated: $0.00145$, $0.00405$ and $0.0013$ for the real part, imaginary one and for the absolute function $|h|$, respectively. The plots of the obtained solution is presented in Figure~\ref{fig:schrodinger}(a) and (b). 

Then, the first $5$ derivatives have been calculated at $t = 0$. The obtained Taylor approximations have been calculated for $100$ random points $x$ at $t = 0.01,$ $0.02,...,$ $0.05$. The obtained $500$ values for real and imaginary parts of $h$ have been added as additional conditions to the model, which was re-initialized, compiled and trained with the same parameters. The final obtained loss has become $9.55e-06$ (a bit bigger that for the standard model due to many additional conditions). The obtained $L_2$ norm at $t = 0.1021$ has become better: $0.000696,$ $0.00356$ and $0.000641$ for the real part, imaginary one and for the absolute function $|h|$, respectively. The plot of the obtained solution is also presented in Figure~\ref{fig:schrodinger}(c) and (d). Again, the proposed scheme allowed to obtain a better solution within the same number of iterations ($\approx31000$ iterations for both the models).

\section{Conclusion}
A new higher order differentiation method has been proposed for time-dependent PDEs. The proposed technique allows to obtain the derivatives with a high precision allowing so to obtain a lot of additional data that can be used for data-driven solution by PINNs. This additional information can be used in different ways, e.g., defining the additional loss functions directly on the derivatives or by adding a high precision Taylor expansion to generate data points inside the domain. 

A well-known DeepXDE software has been used for data-driven solution of several real-life problems under Tensorflow background framework. A way to add new additional data points to the PINN model has been presented in this context. The results have shown that the performance of the deep learning models can be improved in this case. 

It should be noted that the proposed applications of the differentiation techniques presented in this paper have only illustrative purposes. It is already clear that the proposed methods allow to obtain much more data on the dynamics of the system. Since, data is the key value for the neural networks and deep learning, this additional data can be used in different ways in order to improve both the performance and accuracy. Other, more sophisticated models, loss functions and PINNs using the proposed differentiation schemes as well as more efficient differentiation method is the subject of the future works. 

\bibliographystyle{plain}
\bibliography{pde_diff.bib}
\end{document}